\documentclass[11pt,a4paper]{article}
\usepackage[T1]{fontenc}
\usepackage{lmodern}
\usepackage{amsmath,amssymb,amsthm,mathtools,needspace}
\usepackage[textwidth=150mm,textheight=220mm,centering]{geometry}
\usepackage[hidelinks]{hyperref}
\makeatletter
\renewcommand\section{\@startsection{section}{1}{\z@}%
  {-3.5ex \@plus -1ex \@minus -.2ex}%
  {2.3ex \@plus .2ex}%
  {\normalfont\large\bfseries}}
\renewcommand\subsection{\@startsection{subsection}{2}{\z@}%
  {-3.25ex \@plus -1ex \@minus -.2ex}%
  {1.5ex \@plus .2ex}%
  {\normalfont\normalsize\bfseries}}
\makeatother

\newtheorem{theorem}{Theorem}[section]
\newtheorem{proposition}[theorem]{Proposition}
\newtheorem{corollary}[theorem]{Corollary}
\newtheorem{lemma}[theorem]{Lemma}
\newtheorem*{maintheorem}{Theorem A}
\newtheorem*{maintheoremb}{Theorem B}
\newtheorem*{maintheoremc}{Theorem C}
\theoremstyle{remark}
\newtheorem*{introremark}{Remark}
\newtheorem{remark}[theorem]{Remark}
\newtheorem{example}[theorem]{Example}
\newcommand{\Sn}{\mathbb S}
\newcommand{\Tn}{\mathbb T}
\newcommand{\tr}{\operatorname{tr}}
\newcommand{\vol}{\operatorname{Vol}}

\title{Common-Source Assembly for\\
      Spherical Minimal Immersions}
\author{Yongsheng Zhang}
\date{}
\hypersetup{pdftitle={Common-Source Assembly for Spherical Minimal Immersions},
      pdfauthor={Yongsheng Zhang}}

\begin{document}
\maketitle

\begin{abstract}
      We develop a matching principle for assembling maps from a common manifold into spherical minimal immersions.
      The equations couple the coefficients to the induced metric and yield examples with nonminimal constituent immersions.
      We establish existence and uniqueness of matching coefficients under suitable hypotheses
          and sharp rigidity for complete harmonic assemblies over $\Sn^2$.
      For a Clifford family, the first eigenspace determines the underlying fibration.
      The principle also produces weighted minimal profiles that generate minimal products in round spheres,
          with examples formed from finitely many Lawson--Osserman blocks of different degrees over a common Hopf source.
\end{abstract}



\section{Introduction}

      We study when maps from one manifold into different spheres can be assembled into a minimal immersion.
      The coefficients determine both the assembled map and its induced metric, so they also
          determine the Laplacian acting on each block.

\paragraph{The coupled matching problem.}
        Let $f_i:M^m\to\Sn(E_i)$, $1\leq i\leq q$, take values in the unit spheres of Euclidean spaces $E_i$.
        For positive coefficients $t_i$ with $\sum_i t_i=1$, we set
    $$
        A_i=f_i^*g_{\Sn(E_i)},
        \qquad
        F_{\mathbf t}=(\sqrt{t_i}\,f_i)_i=(\sqrt{t_1}\,f_1,\ldots,\sqrt{t_q}\,f_q),
        \qquad
        h_{\mathbf t}=F_{\mathbf t}^*g_{\Sn}
            =\sum_i t_iA_i.
    $$
        If the differentials $df_i$ are jointly injective, then $h_{\mathbf t}$ is a Riemannian metric,
            even when some $A_i$ are singular.

      By Takahashi's criterion \cite{Takahashi}, an isometric immersion is minimal if and only if it satisfies
          $\Delta_hF=-mF$.
      Here we seek coefficients satisfying
    \begin{equation}\label{eq:ordinary-matching}
        \Delta_{h_{\mathbf t}}f_i=-m f_i
        \qquad\text{for every }i.
    \end{equation}
      If $f_i$ is an immersion, its minimality is equivalent to $\Delta_{A_i}f_i=-mf_i$;
       in general,  this is not necessary for a solution of \eqref{eq:ordinary-matching}.
      Example \ref{ex:four} combines two nonminimal torus immersions into a minimal embedding.

\paragraph{Weights and auxiliary factors.}
      We can equip the target sphere with area density $W(z)=\prod_i|z_i|^{k_i}$,
          where $z=(z_i)$ and $k_i\geq0$ are weight exponents.
      As in the minimal product framework of \cite{ZhangMMP}, we regard $F_{\mathbf t}$ as a profile in this weighted sphere.
        With $D=m+\sum_i k_i$, we call $(\mathbf t,\mathbf k)$ matched when
    \begin{equation}\label{eq:mat}
        \Delta_{h_{\mathbf t}}f_i
        =-\left(D-\frac{k_i}{t_i}\right)f_i,
        \qquad \text{for every }i.
    \end{equation}
        At $\mathbf k=0$, this is the minimality equation above.
        For integer $k_i$, take auxiliary minimal isometric immersions $u_i:L_i^{k_i}\to\Sn(V_i)$.
        Matching is then exactly the minimality condition for the tensor map
            $\mathcal F=(\sqrt{t_i}\,f_i\otimes_{\mathbb R}u_i)_i$ on $M\times\prod_iL_i$,
            which arises naturally in the minimal product construction.

        We identify the balance equations with the critical-point equations on the coefficient simplex for
    $$
        J_x(\mathbf t)=\sqrt{\det h_{\mathbf t}(x)}\prod_i t_i^{\frac{k_i}{2}}.
    $$
        The function $\log J_x$ is concave. If every $k_i>0$, it is strictly concave and $J_x$ vanishes on the boundary, so $J_x$ has a unique interior maximizer.
        On a homogeneous source, equivariance makes this maximizer independent of $x$,
            while the isotropy condition in Theorem A supplies harmonicity.

    \begin{maintheorem}[Matching on homogeneous spaces]
        Let $M=G/H$ be connected of dimension $m>0$, with $G$ a compact Lie group and $H$ closed,
            and let $f_i:M\to\Sn(E_i)$ be equivariant maps for orthogonal representations of $G$.
        Assume that their joint differential is injective and that the $H$-fixed subspace of each $E_i$ is $\mathbb R f_i(eH)$.
        Every positive parameter vector $\mathbf k$ admits a unique positive coefficient vector $\mathbf t$,
            with $\sum_i t_i=1$, satisfying \eqref{eq:mat} and depending real analytically on $\mathbf k$.
        The resulting $F_{\mathbf t}$ is stationary for $W$-weighted area;
            for integer parameters, $\mathcal F$ is an ordinary spherical minimal isometric immersion.
    \end{maintheorem}

    \begin{introremark}
        When $\mathbf k=0$, positive matching coefficients need not exist or be unique.
    \end{introremark}

      A geometric application comes from the Lawson--Osserman construction and its generalizations in \cite{XuYangZhang}.
      For the identity map and blocks descending through one fixed Hopf submersion,
          ordinary spherical minimality requires all active descending blocks to have the same eigenvalue.
      The positive weight terms in \eqref{eq:mat} allow arbitrarily many distinct degrees to occur together.
      We construct these weighted minimal profiles in \S~\ref{sec:mmp} and use the minimal product framework
          to obtain minimal immersions in ordinary spheres.

\paragraph{Rigidity at a fixed metric.}
      Every competing minimal isometric immersion has coordinate functions
          in the Laplace eigenspace with eigenvalue equal to the source dimension.
      Over $\Sn^2$, we analyze this eigenspace and determine the linear maps that preserve the induced metric
          within the classical do Carmo--Wallach--T\'oth framework \cite{doCarmoWallach,Toth2002,Toth2006}.
      For positive integers $d_i,k_i$,
      we assemble the complete degree-$d_i$ harmonic maps into a weighted minimal profile
          $\Psi=(\sqrt{t_i}\,\phi_{d_i})_i$ from $\Sn^2$ into the weighted sphere.
     We then attach identity maps of the auxiliary spheres $\Sn^{k_i}$.
      The resulting tensor map descends to a minimal embedding $\mathcal F$ of
      $$
          \mathcal E_{\mathbf d,\mathbf k}
          =\bigl(\Sn^2\times\prod_i\Sn^{k_i}\bigr)/\sim,
          \qquad (x,(u_i))\sim(-x,((-1)^{d_i}u_i)).
      $$

    \begin{maintheoremb}[Rigidity and nonrigidity over $\Sn^2$]
        Every spherical minimal isometric immersion of $\mathcal E_{\mathbf d,\mathbf k}$
                with the metric induced by $\mathcal F$
                    factors as $J\circ\mathcal F$,
                         for some linear isometric inclusion $J$, unless
        $$
            q=2,\qquad k_1=k_2=1,\qquad d_1=d_2=d,
            \qquad p(p+1)=2d(d+1)
        $$
        for some even integer $p\geq2$.
        Each exceptional case admits an explicit one-parameter family of full minimal isometric embeddings not of this form.
    \end{maintheoremb}

        For every integer $N\geq6$, Proposition \ref{prop:metric-families} constructs full minimal embeddings
          of one fixed closed orientable $N$-manifold into round spheres of increasing dimensions.
          These induced metrics are pairwise not locally homothetic.
             For each of these metrics, every full spherical minimal isometric immersion is congruent to the corresponding embedding.

\paragraph{Intrinsic determination of a fibration.}
       We use Clifford systems to construct minimal embeddings for which the first eigenspace determines a fibration of the source.
      We state the case of multiplicities $(3,8)$ here. The general result is Theorem \ref{thm:crec}.
      Let $P_0,\ldots,P_3$ be such an OT--FKM Clifford system and $B^{19}=M_+$ its focal manifold; see \cite{FKM}.
      For $u=(u_0,\ldots,u_3)\in\Sn^3\subset\mathbb R^4$, let $P_u=\sum_{\alpha=0}^3u_\alpha P_\alpha$.
      We identify the unit normal bundle of $B$ with
          $N^{22}=B\times\Sn^3$ via the normal vector $P_ux$ and write $\pi(x,u)=u$ with $\pi_\alpha(x,u)=u_\alpha$.
      We combine $x$, $P_ux$ and the complete degree-$k$ harmonic map of $\Sn^3$ composed with $\pi$.
      Here $E_\mu(h)=\ker(-\Delta_h-\mu)$ denotes the scalar eigenspace.

    \begin{maintheoremc}[A Clifford family]
        For every integer $k\geq2$, we can match these blocks to obtain a full minimal embedding
            $\mathcal G_k:N\to\Sn^{47+(k+1)^2}$ whose induced metric $h_k$ satisfies
        $$
            \lambda_1(h_k)=\frac{66}{k(k+2)},
            \qquad E_{\lambda_1(h_k)}(h_k)=\operatorname{span}\{\pi_0,\ldots,\pi_3\}.
        $$
        The fibers of $\pi$ are exactly the common level sets of the first eigenfunctions,
            so every isometry of $h_k$ preserves the fibration.
        The full isometry group is independent of $k$ as a subgroup of $\operatorname{Diff}(N)$.
    \end{maintheoremc}

      We organize the paper as follows.
      In \S~\ref{sec:mat} we prove the matching principle and develop its applications,
          including Lawson--Osserman profiles over a common Hopf source.
      We study rigidity and nonrigidity in \S~\ref{sec:har},
          then determine the Clifford fibration and its isometry group in \S~\ref{sec:cli}.
      Appendix \ref{app:jet} proves the separation lemma used in the rigidity argument.

\section{A variational matching principle}\label{sec:mat}

      We use $\Delta_h=\operatorname{div}_h\nabla$.
      For a sphere-valued map $f:(M,h)\to\Sn(E)$, we have
          $\Delta_hf=\tau_h(f)-\tr_h(f^*g_{\Sn(E)})f$,
          where $\tau_h(f)$ is the tension field tangent to the sphere.

\subsection{The matching criterion}

      We first separate the matching equations into harmonicity and scalar balance.
      At each point of $M$, the balance equations are the critical point equations for a function of the coefficients.
      We use $\Sigma=\{\mathbf t\in\mathbb R^q:t_i>0,\ \sum_i t_i=1\}$,
          whose closure $\overline\Sigma$ allows $t_i\geq0$.

\begin{theorem}[The matching criterion]\label{thm:mat}
      Let $f_i:M^m\to\Sn(E_i)$, $1\leq i\leq q$, $m>0$, and $A_i=f_i^*g_{\Sn(E_i)}$.
      Assume $\bigcap_i\ker(df_i)_x=\{0\}$ at every $x\in M$.
      For real parameters $k_i\geq0$ and $\mathbf t\in\Sigma$, set
          $h_{\mathbf t}=\sum_i t_iA_i$ and $D=m+\sum_i k_i$.
      The matching equations \eqref{eq:mat} hold if and only if every $f_i$ is harmonic
          with respect to $h_{\mathbf t}$ and
\begin{equation}\label{eq:bal}
          \tr_{h_{\mathbf t}}A_i+\frac{k_i}{t_i}=D
          \qquad\text{for every }i.
\end{equation}
      In a fixed basis of $T_xM$, define
          $J_x(\mathbf t)=\sqrt{\det h_{\mathbf t}(x)}\prod_i t_i^{\frac{k_i}{2}}$.
      Its critical points on $\Sigma$ are exactly the solutions of \eqref{eq:bal} at $x$.
      The function $\log J_x$ is concave.
      If every $k_i>0$, it is strictly concave and $J_x$ has a unique maximizer in $\Sigma$.

      If $\mathbf t_0$ is matched for fixed $\mathbf k\geq0$, then $\mathbf t\in\Sigma$ is matched if and only if
\begin{equation}\label{eq:matching-fiber}
          t_i=t_{0i}\quad\text{for }k_i>0,
          \qquad\text{and}\qquad
          \sum_i(t_i-t_{0i})A_i=0\quad\text{on }M.
\end{equation}
      In particular, all matched coefficient vectors induce the same metric.
\end{theorem}

\begin{proof}
      We split \eqref{eq:mat} into its tangential and normal components to obtain harmonicity and \eqref{eq:bal}.

      Fix $x$ and regard $J_x$ as a function on the positive orthant. We have
$$
          2\partial_{t_i}\log J_x
          =\tr_{h_{\mathbf t}}A_i+\frac{k_i}{t_i}.
$$
      On $\Sigma$, we may vary $\mathbf t$ as $\mathbf t+s(e_i-e_j)$ for small $s$.
      At a critical point of $\log J_x$ on $\Sigma$, we therefore have
          $\partial_{t_i}\log J_x=\partial_{t_j}\log J_x$ for all $i,j$.
      Thus $\tr_{h_{\mathbf t}}A_i+\frac{k_i}{t_i}$ has a common value $c$ for every $i$.
      We multiply by $t_i$ and sum to get
$$
          c=\sum_i t_i\left(\tr_{h_{\mathbf t}}A_i+\frac{k_i}{t_i}\right)
          =\tr_{h_{\mathbf t}}h_{\mathbf t}+\sum_i k_i=D.
$$
      Conversely, if \eqref{eq:bal} holds, the derivative vanishes in every direction tangent to $\Sigma$.
      For $\mathbf v=(v_i)$ with $\sum_i v_i=0$, we compute
\begin{equation}\label{eq:hes}
          2\,\operatorname{Hess}_{\mathbf t}(\log J_x)[\mathbf v,\mathbf v]
          =-\tr\!\left[\left(h_{\mathbf t}^{-\frac{1}{2}}
          \Bigl(\sum_i v_iA_i\Bigr)h_{\mathbf t}^{-\frac{1}{2}}\right)^2\right]
          -\sum_i\frac{k_iv_i^2}{t_i^2}.
\end{equation}
      It follows that $\log J_x$ is concave, strictly so if every $k_i>0$.
      In this case, $J_x$ extends continuously to $\overline\Sigma$, is positive in the interior, and vanishes on the boundary.
      Thus $J_x$ has a unique maximizer in $\Sigma$.

      Finally, suppose $\mathbf t_0$ and $\mathbf t$ are matched.
      For each $x$, both coefficient vectors maximize $\log J_x$, which is therefore constant on the segment joining them by concavity.
      We apply \eqref{eq:hes} with $\mathbf v=\mathbf t-\mathbf t_0$ and obtain \eqref{eq:matching-fiber}.
      Conversely, if \eqref{eq:matching-fiber} holds, the metric and all $D-\frac{k_i}{t_i}$ remain unchanged,
          so $\mathbf t$ also satisfies the matching equations.
\end{proof}

\begin{introremark}[Weighted interpretation]
      For $W(z)=\prod_{k_i>0}|z_i|^{k_i}$, an immersion is $W$-minimal if it is stationary
          under compactly supported variations for $\int W\,d\mu$ on $\{W>0\}$.
      Along $F_{\mathbf t}$, we find that $\nabla^{\Sn}\log W$ is normal to the image, with $i$-th component
          $(\frac{k_i}{t_i}-\sum_jk_j)\sqrt{t_i}f_i$.
      The weighted stationarity equation is
          $\Delta_{h_{\mathbf t}}F_{\mathbf t}+mF_{\mathbf t}
          =\nabla^{\Sn}\log W(F_{\mathbf t})$.
      Componentwise, this is exactly \eqref{eq:mat}.
\end{introremark}

\subsection{Tensor realization and descent}

      As in the minimal product framework of \cite{ZhangMMP}, we attach minimal factors to weighted minimal profiles
          and obtain ordinary spherical minimal products.

      Now suppose the $k_i$ are nonnegative integers.
      For auxiliary maps $u_i:L_i^{k_i}\to\Sn(V_i)$, we record the possible signs by
          $S_i=\{s\in\{\pm1\}:u_i(L_i)\cap s u_i(L_i)\ne\varnothing\}$.
      When $k_i=0$, we take $L_i$ to be one point, $V_i=\mathbb R$, and $u_i=1$, so $S_i=\{1\}$.
      For $u_i=\operatorname{id}_{\Sn^{k_i}}$ with $k_i>0$, we have $S_i=\{\pm1\}$.

\begin{proposition}[Tensor realization and equivariant descent]\label{prop:ten}
      Let $f_i:M^m\to\Sn(E_i)$ and $t_i>0$, with $\sum_i t_i=1$, satisfy \eqref{eq:mat}
          for nonnegative integers $k_i$.
      Assume $h=\sum_i t_if_i^*g_{\Sn(E_i)}>0$ and put $D=m+\sum_i k_i$.
      For minimal isometric immersions $u_i:L_i^{k_i}\to\Sn(V_i)$, the map
$$
          \mathcal F:M\times\prod_iL_i\longrightarrow
          \Sn\!\left(\bigoplus_iE_i\otimes_{\mathbb R}V_i\right),\qquad
          \mathcal F(x,(y_i))=
          \bigl(\sqrt{t_i}\,f_i(x)\otimes_{\mathbb R}u_i(y_i)\bigr)_i
$$
      is a minimal immersion of dimension $D$, with metric $h\oplus\bigoplus_i t_i g_{L_i}$.

      Suppose a discrete group $\Lambda$ acts smoothly, freely, and properly discontinuously on $M$.
      Let $\chi_i:\Lambda\to\{\pm1\}$ be homomorphisms and let $\alpha_{i,\gamma}$ be smooth
          $\Lambda$-actions on $L_i$ such that
          $f_i(\gamma x)=\chi_i(\gamma)f_i(x)$ and $u_i(\alpha_{i,\gamma}y)=\chi_i(\gamma)u_i(y)$.
      Then $\mathcal F$ descends to a minimal immersion of the bundle
          $\mathcal E=(M\times\prod_iL_i)/\Lambda\to M/\Lambda$.
      If $M/\Lambda$ and all $L_i$ are compact, the descended map is an embedding if and only if
          every $u_i$ is an embedding and the following condition holds:
          whenever $f_i(x')=s_i f_i(x)$ for all $i$, with $s_i\in S_i$, there exists
          $\gamma\in\Lambda$ such that $x'=\gamma x$ and $s_i=\chi_i(\gamma)$ for every $i$.
\end{proposition}

\begin{proof}
      Using $\langle f_i,df_i\rangle=\langle u_i,du_i\rangle=0$, we obtain the stated metric.
      Since $\Delta_{t_i g_{L_i}}u_i=-\frac{k_i}{t_i}u_i$, we have $\Delta\mathcal F=-D\mathcal F$ by matching.
      Takahashi's criterion then proves minimality.
      The signs cancel under the diagonal action, which is free and properly discontinuous,
          so $\mathcal F$ descends.

      For unit vectors in real spaces, equality of tensor images is equivalent to
          $f_i(x')=s_i f_i(x)$ and $u_i(y_i')=s_i u_i(y_i)$, with $s_i\in S_i$.
      If the stated condition holds and all $u_i$ are injective, the two source points therefore lie in the same diagonal orbit.
      Conversely, fix $x$ and all but one auxiliary point.
      If the descended map is injective, freeness on $M$ then implies that each $u_i$ is injective.
      For any signs satisfying $f_i(x')=s_i f_i(x)$ with $s_i\in S_i$, we choose auxiliary points realizing them
          and use injectivity of the descended map to obtain the required $\gamma$.
      Compactness turns these injective immersions into embeddings.
\end{proof}

\subsection{Applications of matching}

      The maximizer of $J_x$ may depend on $x$.
      Matching requires a common coefficient vector and harmonicity throughout $M$.
      Symmetry and proportional pullback metrics give two general existence results.
      The Fourier example combines nonminimal inputs, and the Hopf construction matches blocks of distinct degrees.

\subsubsection{Homogeneous sources}

      For equivariant maps on a homogeneous source, the balance equations are independent of $x$.
      The following condition on fixed vectors also ensures harmonicity.
      Given $\rho_i:G\to\operatorname O(E_i)$ and $H\subset G$, write
          $E_i^{\rho_i(H)}=\{v\in E_i:\rho_i(a)v=v\text{ for all }a\in H\}$.

\begin{corollary}[Homogeneous matching]
      Let $M=G/H$ be connected of dimension $m>0$, where $G$ is a compact Lie group and $H$ is closed.
      Let $\rho_i:G\to\operatorname O(E_i)$ be orthogonal representations and
          $f_i:M\to\Sn(E_i)$ smooth equivariant maps with injective joint differential.
      At $o=eH$, assume $E_i^{\rho_i(H)}=\mathbb R f_i(o)$ for every $i$.
      Then each real parameter vector $\mathbf k$ with $k_i>0$ admits a unique matching coefficient vector
          $\mathbf t_*(\mathbf k)$, depending real analytically on $\mathbf k$.
\end{corollary}

\begin{proof}
      Since $h_{\mathbf t}$ is $G$-invariant, we have
          $\tau_{h_{\mathbf t}}(f_i)(o)\in E_i^{\rho_i(H)}\cap f_i(o)^\perp=\{0\}$ by equivariance.
      Thus each $f_i$ is harmonic everywhere.
      As $\tr_{h_{\mathbf t}}A_i$ is constant on $M$, we can apply Theorem \ref{thm:mat} at $o$
          to obtain the unique global matching vector.
      Real analytic dependence follows from \eqref{eq:hes} and the analytic implicit function theorem.
\end{proof}

      This proves the existence and uniqueness assertions of Theorem A.

\begin{introremark}
      Take $G=\operatorname{SU}(3)$, $H=\operatorname{S}(\operatorname U(1)\times\operatorname U(2))$,
          and let $G$ act by $\rho(g)A=gAg^*$ on the real space $E$ of traceless Hermitian $3\times3$ matrices
          with the trace inner product.
      We have $E^{\rho(H)}=\mathbb R\operatorname{diag}(2,-1,-1)$.
      The normalized orbit map is the standard embedding $\mathbb{CP}^2\to\Sn^7$ and satisfies the hypothesis above.
\end{introremark}

\subsubsection{Homothetic eigenmaps}

      When all pullback metrics are proportional to one metric, we can determine the matching coefficients
          from a single scalar equation.
      Examples come from complete Laplace eigenspaces on spheres and projective spaces with their standard metrics
          \cite[\S~3]{Takahashi}.
      By symmetry, the evaluation maps have constant length and pullback metrics proportional to the source metric.

\begin{corollary}[Matching homothetic eigenmaps]\label{cor:eig}
      Let $(M^m,g)$, $m>0$, admit sphere-valued maps $f_i$ satisfying
$$
          f_i^*g_{\Sn(E_i)}=\frac{\lambda_i}{m}g,\qquad
          \Delta_g f_i=-\lambda_i f_i,
$$
      where $\lambda_i\geq0$ and at least one $\lambda_i>0$.
      For real $k_i>0$, set $D=m+\sum_i k_i$.
      There is a unique $\sigma>\max_i\lambda_i$ such that
\begin{equation}\label{eq:eig}
          \sum_i\frac{k_i\lambda_i}{\sigma-\lambda_i}=m.
\end{equation}
      The unique matching coefficients and metric are
          $t_i=\frac{k_i\sigma}{D(\sigma-\lambda_i)}$ and $h_{\mathbf t}=\frac{\sigma}{D}g$.
\end{corollary}

\begin{proof}
      We have $h_{\mathbf t}=cg$.
      From the pullback identities, we obtain $c=m^{-1}\sum_i t_i\lambda_i>0$.
      We set $\sigma=Dc$ and rewrite the matching equations as
          $t_i=\frac{k_i\sigma}{D(\sigma-\lambda_i)}$.
      Since the $t_i$ are positive, we must have $\sigma>\max_i\lambda_i$.
      Using $\sum_i t_i=1$, we reduce the problem to \eqref{eq:eig}.
      The left-hand side decreases strictly from $+\infty$ to zero on $(\max_i\lambda_i,\infty)$,
          so the solution is unique.
      Conversely, choose this root and define $t_i$ by the formula above.
      We then have $t_i>0$, $\sum_i t_i=1$ and $m^{-1}\sum_i t_i\lambda_i=\frac{\sigma}{D}$.
      Hence these coefficients solve the matching equations.
\end{proof}

      For each $\lambda_i>0$, we have $\Delta_{\frac{\lambda_i}{m}g}f_i=-mf_i$,
          so the input itself is minimal for its induced metric.
      For example, the degree-one and degree-two harmonic maps of the unit sphere $(\Sn^2,g_0)$ have eigenvalues $2$ and $6$.
      If we choose $(k_1,k_2)=(2,1)$, we obtain $\sigma=10$, $t_1=t_2=\frac12$, and $h_{\mathbf t}=2g_0$.
      With $k_1=k_2=0$, no positive combination of these two maps is minimal in a round sphere,
          since their eigenvalues for any $cg_0$ remain different.

\subsubsection{Fourier maps}

      Fourier constructions on flat tori are classical; see \cite{Bryant,LvWangXie2024,ParkUrakawa}.
      We apply Theorem \ref{thm:mat} to combine two nonminimal immersions into a minimal embedding.

\begin{example}[Two nonminimal inputs]\label{ex:four}
      We take $k_1=k_2=0$, $t_1=t_2=\frac12$, even $j\geq2$, and $0<\alpha<1$ with $\alpha\ne\frac12$.
      On $\Tn^2=(\mathbb R/2\pi\mathbb Z)^2$, define $f,g:\Tn^2\to\Sn^3$ by
$$
          \begin{aligned}
          f(x,y)&=\bigl(\sqrt\alpha\,e^{\mathrm i(x+jy)},
          \sqrt{1-\alpha}\,e^{\mathrm i(x-jy)}\bigr),\\
          g(x,y)&=\bigl(\sqrt{1-\alpha}\,e^{\mathrm i(jx+y)},
          \sqrt\alpha\,e^{\mathrm i(jx-y)}\bigr).
          \end{aligned}
$$
      Both inputs cover product tori in $\Sn^3$ with constant nonzero mean curvature.
      With $c=j(2\alpha-1)$, their induced metrics are
$$
          A=\begin{pmatrix}1&c\\c&j^2\end{pmatrix},\qquad
          B=\begin{pmatrix}j^2&-c\\-c&1\end{pmatrix}.
$$
      Both are positive definite, and $F=\frac{1}{\sqrt2}(f,g)$ induces $h=\frac{A+B}{2}=\frac{(j^2+1)I}{2}$.
      Since $\Delta_h=\frac{2(\partial_x^2+\partial_y^2)}{j^2+1}$, we get $\Delta_hf=-2f$ and $\Delta_hg=-2g$,
          so $F$ is minimal.
      For either input, however, the coordinate eigenvalues for its own induced metric are
          $\frac{1}{\alpha}$ and $\frac{1}{1-\alpha}$.
      Since $\alpha\ne\frac12$, neither immersion is minimal.

      To check injectivity, suppose $F(x,y)=F(x',y')$ and put $\delta x=x'-x$, $\delta y=y'-y$.
      We have $\delta x\pm j\delta y\in2\pi\mathbb Z$ and $j\delta x\pm\delta y\in2\pi\mathbb Z$.
      Adding and subtracting these relations, we obtain $\delta x,\delta y\in\pi\mathbb Z$.
      As $j$ is even, the original relations now imply $\delta x,\delta y\in2\pi\mathbb Z$, so $F$ is injective.
      The four frequencies are distinct up to sign, so their real and imaginary coordinate functions are linearly independent.
      Thus $F:\Tn^2\to\Sn^7$ is a full minimal embedding.
\end{example}

\subsubsection{Matching over a common Hopf source}\label{sec:mmp}

      In joint work with Xu and Yang \cite[Theorems~2.5, 2.8 and 2.10]{XuYangZhang}, we characterize nontrivial Lawson--Osserman maps
          with equal nonzero singular values as Hopf pullbacks of homothetic minimal immersions of
          complex or quaternionic projective spaces, or of the octonionic projective line.
      These immersions may be chosen arbitrarily from the corresponding classical moduli spaces.
      We combine the identity map with finitely many such pullbacks through one fixed Hopf submersion.
      Their pullback metrics vanish in the fiber directions, so they are not proportional to the metric of the identity block.

      Let $L_\alpha$ be the pullback eigenvalues and $B$ the horizontal coefficient of the induced metric.
      In an ordinary round sphere, matching requires $\frac{L_\alpha}{B}=n$ for every active pullback block,
          where $n$ is the dimension of the Hopf source.
      Thus all active pullback blocks in this construction must have the same eigenvalue, and hence the same degree.
      This obstruction was observed by the author and Ling Yang during their visit to MSRI in 2017.
      With positive weight exponents, we instead have
$$
          \frac{L_\alpha}{B}+\frac{k_\alpha}{t_\alpha}=D,
          \qquad D=n+\sum_{\alpha=0}^s k_\alpha.
$$
      The additional weight terms allow us to match any finite collection of these blocks.

\begin{theorem}[Multiple Lawson--Osserman blocks]\label{thm:lo-products}
      Let $\pi:\Sn^n\to(P^r,g_P)$, $0<r<n$, be a standard Hopf submersion.
      Its fibers are totally geodesic and $\pi(-x)=\pi(x)$.
	      Write $g_{\Sn^n}=g_V+g_H$, where $g_H=\pi^*g_P$.
      For $1\leq\alpha\leq s$, let $\bar f_\alpha:P\to\Sn^{N_\alpha}$ be linearly full maps satisfying
      $$
          \bar f_\alpha^*g_{\Sn^{N_\alpha}}=\frac{L_\alpha}{r}g_P,
          \qquad
          \Delta_{g_P}\bar f_\alpha=-L_\alpha\bar f_\alpha,
      $$
      and set $f_\alpha=\bar f_\alpha\circ\pi$.
      Fix positive integers $k_0,k_1,\ldots,k_s$.
      For $t_\alpha>0$, $\sum_{\alpha=0}^s t_\alpha=1$, define
      $$
          B_{\mathbf t}=t_0+\sum_{\alpha=1}^s\frac{L_\alpha}{r}t_\alpha,
          \qquad
          J(\mathbf t)=t_0^{\frac{n-r+k_0}{2}}B_{\mathbf t}^{\frac{r}{2}}
              \prod_{\alpha=1}^s t_\alpha^{\frac{k_\alpha}{2}}.
      $$
      Then $J$ has a unique maximizer $\mathbf t_*=(t_0,\ldots,t_s)$
          in the open coefficient simplex.

      The profile $F_{\mathbf t_*}(x)=(\sqrt{t_0}\,x,(\sqrt{t_\alpha}\,f_\alpha(x))_{\alpha=1}^s)$
          is minimal in the weighted sphere with density
          $W(z)=\prod_{\alpha=0}^s|z_\alpha|^{k_\alpha}$.

      With $u_\alpha\in\Sn^{k_\alpha}$, the map
      $$
          \mathcal L_{\mathbf t_*}(x,u_0,\ldots,u_s)
          =\left(\sqrt{t_0}\,x\otimes u_0,
              \bigl(\sqrt{t_\alpha}\,f_\alpha(x)\otimes u_\alpha\bigr)_{\alpha=1}^s\right)
      $$
      descends to a linearly full minimal embedding
      $$
          \frac{\Sn^n\times\prod_{\alpha=0}^s\Sn^{k_\alpha}}
          {(x,u_0,u_1,\ldots,u_s)\sim(-x,-u_0,u_1,\ldots,u_s)}
          \longrightarrow\Sn^M,
      $$
      where $M=(n+1)(k_0+1)+\sum_{\alpha=1}^s(N_\alpha+1)(k_\alpha+1)-1$.
      Its covering metric is
      $$
          t_0g_V+B_{\mathbf t_*}g_H
          \oplus\bigoplus_{\alpha=0}^s t_\alpha g_{\Sn^{k_\alpha}}.
      $$
\end{theorem}

\begin{proof}
      The weighted volume density relative to the round metric is $J$, independent of the source point.
      By Theorem \ref{thm:mat}, it has a unique interior maximizer.
      The Hopf fibers are totally geodesic, so for $h_{\mathbf t}=t_0g_V+B_{\mathbf t}g_H$ we have
$$
          \Delta_{h_{\mathbf t}}x
          =-\left(\frac{n-r}{t_0}+\frac r{B_{\mathbf t}}\right)x,
          \qquad
          \Delta_{h_{\mathbf t}}f_\alpha=-\frac{L_\alpha}{B_{\mathbf t}}f_\alpha.
$$
      Every block is therefore harmonic for $h_{\mathbf t}$.
      Applying Theorem \ref{thm:mat}, we obtain the unique matched vector $\mathbf t_*$
          and hence a weighted minimal profile.
      We obtain the tensor immersion and its metric from Proposition \ref{prop:ten}.

      In the descent criterion of that proposition, the identity block forces $x'=\varepsilon x$,
          $\varepsilon\in\{\pm1\}$.
      Since $\pi(-x)=\pi(x)$, every other block is even and its sign must be $+1$.
      We thus obtain exactly the stated quotient and an embedding.
      We separate the blocks by changing one auxiliary sign at a time.
      Since each block is full, so is the assembly.
\end{proof}

      The proof does not require equivariance of the maps $\bar f_\alpha$.
      For fixed $L_\alpha$ and $k_\alpha$, the matching coefficients and induced metric remain unchanged
          as these maps vary within their moduli spaces of minimal immersions.

      For example, take the Hopf submersion from $\Sn^3$ onto the round $2$-sphere of radius $\frac12$,
          and the complete degree-one and degree-two harmonic maps on the base.
      Their pullbacks have $L_1=8$ and $L_2=24$.
      Together with the identity block, the choices
$$
          (k_0,k_1,k_2)=(8,1,2),\qquad
          (t_0,t_1,t_2)=\left(\frac23,\frac1{12},\frac14\right)
$$
      give $B_{\mathbf t}=4$ and $D=14$.
      The block eigenvalues $(2,2,6)$ and weight contributions $(12,12,8)$ sum to $14$ in each block.
      Thus the spherical degrees $1,2,4$ occur simultaneously in a matched profile.

\begin{remark}[The two-block coefficient]
      In Theorem \ref{thm:lo-products}, take $s=1$ and $L=L_1>n$, and let $t_*$ be the identity coefficient.
      For the ordinary Lawson--Osserman construction \cite{LawsonOsserman,XuYangZhang}, this coefficient is
          $t_{\mathrm{LO}}=\frac{L(n-r)}{n(L-r)}$.
      The two coefficients agree precisely when $k_0:k_1=L(n-r):r(L-n)$.
      Indeed, at $t_{\mathrm{LO}}$ the unweighted density is stationary, so
          $\partial_t\log J=\frac{k_0}{2t_{\mathrm{LO}}}-\frac{k_1}{2(1-t_{\mathrm{LO}})}$.
      This vanishes precisely for the stated ratio.
      By uniqueness, the critical point is $t_*$.
\end{remark}

\section{Rigidity and nonrigidity for harmonic assemblies}\label{sec:har}

      We apply the matching principle to complete spaces of spherical harmonics.
      We construct these embeddings, compute the spectra of their induced metrics, and ask whether the metrics determine the embeddings.
      Over $\Sn^2$, we obtain a sharp rigidity criterion and construct nonrigid realizations in the exceptional cases.

\subsection{Harmonic blocks and quotient spectrum}

      Let $g_0$ be the unit round metric on $\Sn^m$, $m\geq1$.
      For $d\geq1$, we use the full space $E_d$ of real homogeneous harmonic polynomials of degree $d$ on $\mathbb R^{m+1}$, restricted to $\Sn^m$.
      Using an $L^2(g_0)$-orthonormal basis $(Y_\alpha)$, define
          $\phi_d=\sqrt{\frac{\vol(\Sn^m,g_0)}{\dim E_d}}\,(Y_\alpha)_\alpha$.
      Rotations act orthogonally on $E_d$, and $\phi_d$ is equivariant for this action.
      By transitivity, $|\phi_d|$ is constant, and the normalization gives $|\phi_d|=1$.
      The pullback metric is also rotation invariant.
      Since the isotropy group acts irreducibly on each tangent space, we have $\phi_d^*g=cg_0$.
      From $\Delta_{g_0}\phi_d=-d(d+m-1)\phi_d$ and $|\phi_d|=1$, we get
          $mc=\lambda_d=d(d+m-1)$ and hence
\begin{equation}\label{eq:har}
          \phi_d^*g_{\Sn(E_d)}=\frac{\lambda_d}{m}g_0,
          \qquad \Delta_{g_0}\phi_d=-\lambda_d\phi_d.
\end{equation}
      The standard map $\phi_d:\Sn^m\to\Sn(E_d)$ therefore satisfies Corollary \ref{cor:eig}; see \cite[\S~1]{Toth2002} for the classical construction.
      We use the matching coefficients in the tensor construction of Proposition \ref{prop:ten}.
      The spectrum below identifies the possible coordinate functions of any competing minimal immersion.

\begin{proposition}[Complete harmonic assembly and quotient]\label{prop:har}
      Fix $d_1,\dots,d_q\geq1$ and positive integers $k_1,\dots,k_q$.
      Set $\lambda_i=d_i(d_i+m-1)$ and $D=m+\sum_i k_i$, and take $\sigma,t_i$ from Corollary \ref{cor:eig}.
      The profile $\Psi=(\sqrt{t_i}\,\phi_{d_i})_i$ is minimal in the weighted sphere
          with density $W(z)=\prod_i|z_i|^{k_i}$ and induces $h=\frac{\sigma}{D}g_0$.

      Assume $m\geq2$.
      Then $\Psi$ embeds $\Sn^m$ if some $d_i$ is odd, and descends to an embedding of $\mathbb RP^m$ otherwise.
      With identity factors $u_i:\Sn^{k_i}\to\Sn^{k_i}$, the tensor map
\begin{equation}\label{eq:hlift}
          \mathcal F(x,u_1,\dots,u_q)
          =\bigl(\sqrt{t_i}\,\phi_{d_i}(x)\otimes u_i\bigr)_{i=1}^q
\end{equation}
      descends to a linearly full minimal embedding of
\begin{equation}\label{eq:hquot}
          \mathcal E_{\mathbf d,\mathbf k}:=
          \bigl(\Sn^m\times\prod_i\Sn^{k_i}\bigr)/\langle\iota\rangle,
          \qquad
          \iota(x,u_1,\dots,u_q)
          =(-x,(-1)^{d_1}u_1,\dots,(-1)^{d_q}u_q),
\end{equation}
      with covering metric $h\oplus\bigoplus_i t_i g_{\Sn^{k_i}}$.

      A complete eigenbasis on this quotient consists of products of spherical harmonics.
      Write $p$ for the degree on the base $\Sn^m$ and $\ell_i$ for the degree on the $i$-th auxiliary sphere $\Sn^{k_i}$.
      The eigenvalues of $-\Delta$ are
\begin{equation}\label{eq:hspec}
          \frac{\sigma}{D}\nu_{p,\boldsymbol\ell}
          =p(p+m-1)+\sum_i\frac{\sigma-\lambda_i}{k_i}
          \ell_i(\ell_i+k_i-1),
          \qquad p+\sum_i d_i\ell_i\equiv0\pmod2,
\end{equation}
      where $p,\ell_1,\ldots,\ell_q$ are nonnegative integers.
\end{proposition}

\begin{proof}
      Corollary \ref{cor:eig} gives the matching coefficients and metric.
      If $m\geq2$ and $y\ne\pm x$, choose orthonormal $a,b\in x^\perp$ with $\langle a,y\rangle\ne0$ and $\langle b,y\rangle=0$.
      Then $P(z)=\operatorname{Re}(\langle a,z\rangle+\mathrm i\langle b,z\rangle)^d$ is harmonic,
          with $P(x)=0\ne P(y)$.
      Hence $\phi_d(y)\ne\pm\phi_d(x)$.
      Together with $\phi_d(-x)=(-1)^d\phi_d(x)$, this identifies the fibers of $\Psi$.
      The induced metric is positive, so compactness gives the asserted embeddings.

      These signed equalities verify the embedding criterion in Proposition \ref{prop:ten}, so the tensor map embeds the quotient \eqref{eq:hquot}.
      Changing $u_i$ to $-u_i$ isolates the $i$-th tensor block, and each factor spans its ambient space, so $\mathcal F$ is full.
      Products of spherical harmonics form a complete eigenbasis on the covering product.
      After the metric rescaling, their eigenvalues add as in \eqref{eq:hspec}. Invariance under $\iota$ is exactly the stated parity condition.
\end{proof}

\begin{remark}[Topology of the quotient]\label{rem:htop}
      The space \eqref{eq:hquot} is a bundle over $\mathbb RP^m$ with fiber $\prod_i\Sn^{k_i}$.
      A generating loop in the base acts antipodally on the $i$-th factor when $d_i$ is odd.
      Since the antipodal map on $\Sn^n$ has degree $(-1)^{n+1}$,
          the quotient is orientable precisely when $m+1+\sum_i d_i(k_i+1)\equiv0\pmod2$.
      If $m\geq2$ and every $k_i\geq2$, the covering product is simply connected and $\pi_1(\mathcal E_{\mathbf d,\mathbf k})=\mathbb Z_2$.
\end{remark}

\begin{introremark}
      The auxiliary factors need not be round spheres.
      If $k_i=2$ and $d_i$ is even, we may replace the $i$-th sphere by a Lawson minimal embedding
          $u_i:\Sigma\to\Sn^3$ of any prescribed genus \cite[Theorem~2]{Lawson1970}.
      Proposition \ref{prop:ten} gives a minimal embedding with the same matching coefficients,
          taking $\iota$ to act trivially on $\Sigma$.
      Formula \eqref{eq:hspec} and the rigidity analysis below concern the round auxiliary factors.
\end{introremark}
\subsection[Rigidity and nonrigidity over S2]{Rigidity and nonrigidity over $\Sn^2$}

      Let $H_d=E_d$ denote the degree-$d$ harmonics on $\Sn^2$, with evaluation map $\phi_d$ from \eqref{eq:har}.
      The following separation lemma determines the coefficient maps in the rigidity proof.
      The case of equal degrees and symmetric operators is classical; see \cite{Calabi1967,Toth2002}.
      We prove the statement for arbitrary operators in Appendix \ref{app:jet}.

\begin{lemma}[First-jet separation]\label{lem:jet}
      For $p,q\geq1$, let $T:H_p\to H_q$ be linear.
      If $\langle T\phi_p,\phi_q\rangle\equiv0$ and $\langle T\,d\phi_p,\phi_q\rangle\equiv0$, then $T=0$.
\end{lemma}

\begin{theorem}[Rigidity for complete harmonic assemblies over $\Sn^2$]\label{thm:rig}
      In Proposition \ref{prop:har}, take $m=2$, $q\geq1$, and $d_i,k_i\geq1$.
      Assume that we are not in the case $q=2$, $k_1=k_2=1$, $d_1=d_2=d$ with
          $p(p+1)=2d(d+1)$ for some even integer $p\geq2$.
      Set $Q=\sum_i(2d_i+1)(k_i+1)$, and let $h$ be the metric induced by $\mathcal F$ on $\mathcal E_{\mathbf d,\mathbf k}$.
      Every minimal isometric immersion $G:(\mathcal E_{\mathbf d,\mathbf k},h)\to\Sn^M$ equals $J\circ\mathcal F$,
          where $\mathcal F$ is the tensor embedding \eqref{eq:hlift}
          and $J:\mathbb R^Q\to\mathbb R^{M+1}$ is a linear isometric inclusion.
\end{theorem}

      Proposition \ref{prop:mix} gives nonrigid realizations in the excluded cases.

\begin{proof}
      \emph{Step 1: the eigenspace at eigenvalue $D$.}
      Set $D=2+\sum_i k_i$, $\lambda_i=d_i(d_i+1)$, and $\delta_i=\sigma-\lambda_i$,
          with $\sigma$ from \eqref{eq:eig}.
      In \eqref{eq:hspec}, $p$ and $\ell_i$ are the degrees on the base and the $i$-th auxiliary sphere.
      At eigenvalue $D$, its right-hand side equals $\sigma$.
      From \eqref{eq:eig}, we have $\delta_i\geq\frac{k_i\lambda_i}{2}$ and hence
$$
          \frac{2(k_i+1)}{k_i}\delta_i-\sigma
          =\frac{k_i+2}{k_i}\delta_i-\lambda_i\geq\frac{k_i\lambda_i}{2}>0,
$$
      so a $D$-eigenmode has $\ell_i\in\{0,1\}$.
      For $i\ne j$, we also obtain
$$
          4\geq\sigma\left(\frac1{\delta_i}+\frac1{\delta_j}\right)
          \geq\frac{4\sigma}{\delta_i+\delta_j},
$$
      hence $\delta_i+\delta_j\geq\sigma$.
Two nonzero $\ell_i$ force equality in both inequalities:
          $q=2$, $k_1=k_2=1$, and $\delta_1=\delta_2$, hence $d_1=d_2$.
      The modes with at most one nonzero auxiliary degree are
$$
          \text{(1)}\quad(p,\boldsymbol\ell)=(d_i,e_i),
          \qquad
          \text{(2)}\quad\boldsymbol\ell=0,\quad p(p+1)=\sigma,\quad p\text{ even},
$$
      where $e_i$ is the $i$-th standard basis vector of $\mathbb R^q$.
      The case with two circles above also admits $\text{(3)}\ (p,\boldsymbol\ell)=(0,(1,1))$.

      \emph{Step 2: eliminating the additional terms.}
      By Step 1, the lift of $G$ has the form
$$
          \widetilde G=\sum_iB_i(\phi_{d_i}\otimes u_i)+C(x)+Z(\mathbf u),
$$
      where $B_i:H_{d_i}\otimes\mathbb R^{k_i+1}\to\mathbb R^{M+1}$ are linear,
          $C$ occurs only in (2), and $Z=R(u_1\otimes u_2)$ only in (3), with $R$ linear.
      Isometry gives $\sum_i\frac{|d^{(i)}\widetilde G|^2}{k_i}=1=|\widetilde G|^2$, using unit round metrics for $d^{(i)}$.
We average over the auxiliary spheres with normalized volume, eliminating cross terms.
      In the averaged energy, the squared norm of each term in (1) has coefficient one, that of $Z$ two, and that of $C$ zero.
      We obtain $\int|Z|^2\,d\mathbf u=|C(x)|^2$.
      Our hypothesis rules out the simultaneous occurrence of (2) and (3), so the last equality implies $C=Z=0$.

      \emph{Step 3: the linear maps preserve inner products.}
      Put $B_{i,w}\xi=B_i(\xi\otimes w)$.
      For $u_i\perp w$, isometry implies $|B_{i,w}\phi_{d_i}|^2=t_i|w|^2$.
      Polarization and mixed vertical isometry imply
$$
          \langle B_{i,w}\phi_{d_i},B_{j,z}\phi_{d_j}\rangle
          =\begin{cases}t_i\langle w,z\rangle,&i=j,\\0,&i\ne j.\end{cases}
$$
      We also have $\langle B_{i,w}\,d\phi_{d_i},B_{j,z}\phi_{d_j}\rangle=0$.
      It suffices to verify this for $w\ne0$.
      For $i=j$, use the mixed base--sphere metric at $u_i=\frac{w}{|w|}$, $z\perp w$, averaging over the other signs.
      For $z=w$, differentiate the norm identity, then extend by linearity.
      For $i\ne j$, take $u_i=\frac{w}{|w|}$, $u_j\perp z$ in the mixed base--sphere identity and subtract its value at $-u_i$.
      We apply Lemma \ref{lem:jet} to $B_{i,z}^*B_{i,w}-t_i\langle z,w\rangle I$ for $i=j$,
          and to $B_{j,z}^*B_{i,w}$ otherwise.
      We obtain $B_i^*B_i=t_iI$ and $B_i^*B_j=0$ for $i\ne j$.
      Thus $J((\xi_i)_i)=\sum_i t_i^{-\frac{1}{2}}B_i\xi_i$ is a linear isometric inclusion.
      On the quotient, we obtain $G=J\circ\mathcal F$. This completes the proof.
\end{proof}

\begin{proposition}[Nonrigid realizations]\label{prop:mix}
      Let $q=2$, $k_1=k_2=1$, $d_1=d_2=d$, and suppose an even integer $p\geq2$ satisfies $p(p+1)=2d(d+1)$.
      For $0<\varepsilon<\frac12$, the map
$$
          G_\varepsilon(x,u,v)=\bigl(
          \sqrt{\tfrac12-\varepsilon}\,\phi_d(x)\otimes u,
          \sqrt{\tfrac12-\varepsilon}\,\phi_d(x)\otimes v,
          \sqrt\varepsilon\,u\otimes v,
          \sqrt\varepsilon\,\phi_p(x)\bigr)
$$
      descends to a full minimal isometric embedding of $(\mathcal E_{(d,d),(1,1)},h)$ into $\Sn^{8d+2p+8}$,
          where $h$ is the metric induced by the canonical embedding.
      None of these embeddings factors as $J\circ\mathcal F$ for a linear isometric inclusion $J$.
      For a fixed source parametrization, the embeddings for distinct $\varepsilon$ are inequivalent under orthogonal transformations of the target.
\end{proposition}

\begin{proof}
      Set $\lambda=d(d+1)$.
      On the cover, we have $|G_\varepsilon|=1$ and
$$
          G_\varepsilon^*g_{\Sn}
          =\left((1-2\varepsilon)\frac\lambda2+\varepsilon\lambda\right)g_0
          +\frac12(g_{\Sn^1}+g_{\Sn^1})
          =\widetilde h.
$$
      Since $\Delta_{\widetilde h}G_\varepsilon=-4G_\varepsilon$, minimality follows from Takahashi's criterion.
      Since $p$ is even, the deck involution $(x,u,v)\mapsto(-x,(-1)^du,(-1)^dv)$ preserves every block.
The map $G_\varepsilon$ descends to the quotient, where its first two blocks are $\sqrt{1-2\varepsilon}$ times the canonical embedding and ensure embeddedness.
      By changing the signs of $u,v$ independently, we see that the four coordinate spaces are linearly independent.
      Each block is full, so $G_\varepsilon$ is full.
      Any factorization $G_\varepsilon=J\circ\mathcal F$ would place its image in a subspace of dimension at most
          $4(2d+1)<8d+2p+9$, contradicting fullness.

      Let $\Pi$ be the scalar $L^2$-projection onto the span of the coordinate functions of $u\otimes v$,
          using probability volume on the quotient.
      We have $\sum_\alpha\|\Pi(G_\varepsilon)_\alpha\|_{L^2}^2=\varepsilon$.
      Orthogonal transformations of the target preserve this sum, so embeddings with distinct parameters are inequivalent.
\end{proof}

\paragraph{Families on fixed manifolds.}
      We now vary one harmonic degree while keeping the quotient manifold fixed.
      For an $\Sn^2$ base, we obtain infinitely many metrics with rigid minimal realizations in every total dimension at least six.

\begin{proposition}[Metric families and rigidity]\label{prop:metric-families}
      Fix integers $m,a,b\geq2$, and set $N=m+a+b$ and $L_r=2r(2r+m-1)$ for $r\geq1$.
      Let $\sigma_r,t_{1,r},t_{2,r}$ be given by Corollary \ref{cor:eig} for
          $(\lambda_1,\lambda_2)=(m,L_r)$ and $(k_1,k_2)=(a,b)$.
      The assemblies of degrees $(1,2r)$ give full minimal embeddings $F_r$ of the fixed manifold
$$
          \mathcal E^N_{m,a,b}
          =\frac{\Sn^m\times\Sn^a\times\Sn^b}
          {(x,u,v)\sim(-x,-u,v)}
$$
      into round spheres.
      Their induced metrics $h_r$ have covering metrics
$$
          \widetilde h_r=\frac{\sigma_r}{N}g_{\Sn^m}
          +t_{1,r}g_{\Sn^a}+t_{2,r}g_{\Sn^b},
$$
      and are pairwise not locally homothetic.

      When $m=2$, every minimal isometric immersion
          $G:(\mathcal E^N_{2,a,b},h_r)\to\Sn^M$ equals $J\circ F_r$ for a linear isometric inclusion $J$.
      Taking $a=2$ and $b=N-4$, we obtain, for every $N\geq6$, a fixed closed orientable
          $N$-manifold with fundamental group $\mathbb Z_2$ and infinitely many such metrics, each admitting a rigid full minimal embedding
          $F_r:\mathcal E^N_{2,2,N-4}\to\Sn^{8+(4r+1)(N-3)}$.
\end{proposition}

\begin{proof}
      By Proposition \ref{prop:har}, we obtain the embeddings and covering metrics.
      Each summand in \eqref{eq:eig} is less than $m$, so the squared radius ratios satisfy
          $R_{1,r}=\frac{\sigma_r-m}{a}>1$ and $R_{2,r}=\frac{\sigma_r-L_r}{b}>\frac{L_r}{m}>1$.
      Thus the base sphere has the largest radius.
      The other two factors are distinguished by dimension if $a\ne b$, and by $t_{1,r}<t_{2,r}$ if $a=b$.
      The local de Rham decomposition therefore makes $R_{1,r}$ invariant under local homothety.
      The left side of \eqref{eq:eig} decreases in $\sigma$ and increases in $L_r$,
          so $\sigma_r$ and $R_{1,r}$ strictly increase with $r$.
      Thus the metrics are pairwise not locally homothetic.

      For $m=2$, the degrees $1$ and $2r$ are distinct, so rigidity follows from Theorem \ref{thm:rig}.
      With $a=2$ and $b=N-4$, we obtain orientability and fundamental group $\mathbb Z_2$ from Remark \ref{rem:htop}.
      The target dimension is $3(a+1)+(4r+1)(b+1)-1=8+(4r+1)(N-3)$.
\end{proof}

\subsection{Nonrigidity for higher-dimensional bases}

      The rigidity proof over $\Sn^2$ uses Lemma \ref{lem:jet} to rule out nonzero quadratic relations among the coordinates.
      Such relations can occur on higher-dimensional bases.
      We use them to deform repeated blocks while keeping the induced metric and target sphere fixed.
This construction requires no additional eigenfunctions, unlike Proposition \ref{prop:mix}.

\begin{proposition}[Deformations of repeated blocks]
      Let $f:M\to\Sn(V)$ be smooth and full, $q\geq2$, $k_i\geq1$, and $t_i>0$ with $\sum_i t_i=1$.
      Suppose $F=(\sqrt{t_i}\,f\otimes u_i)_i$, with $u_i\in\Sn^{k_i}$,
          is a minimal immersion of $M\times\prod_i\Sn^{k_i}$ with induced metric $g$.
      Write $\mathcal K_0(f)=\{R\in\operatorname{End}(V):R^*=R,\ \langle Rf,f\rangle\equiv0\}$.
      If $R_i\in\mathcal K_0(f)$, $\sum_i t_iR_i=0$, and $I+R_i>0$, then
$$
          F_{\mathbf R}(x,\mathbf u)
          =\bigl(\sqrt{t_i}\,(I+R_i)^{\frac{1}{2}}f(x)\otimes u_i\bigr)_{i=1}^q
$$
      is a full minimal isometric immersion for $g$, and an embedding whenever $F$ is.
      For a fixed source parametrization, distinct tuples yield immersions that are inequivalent under orthogonal transformations of the target.
      Near zero, these tuples form a family of dimension $(q-1)\dim\mathcal K_0(f)$.
\end{proposition}

\begin{proof}
      We have $F_{\mathbf R}=A_{\mathbf R}F$, where $A_{\mathbf R}=\bigoplus_i((I+R_i)^{\frac{1}{2}}\otimes I)$ is constant and invertible.
      Since $\langle R_if,f\rangle=0$, the norm and the metric along each auxiliary sphere are unchanged.
Differentiating $\langle R_if,f\rangle=0$, we find that the mixed terms between the base and auxiliary spheres still vanish.
      The change in the base metric is $\sum_i t_i\langle R_i\,df(X),df(Y)\rangle=0$.
      Thus $F_{\mathbf R}^*g_{\Sn}=g$ and $\Delta_gF_{\mathbf R}=-(\dim M+\sum_i k_i)F_{\mathbf R}$.
      It follows that $F_{\mathbf R}$ is minimal.
      Since $f$ is full, varying the auxiliary factors shows that $F$ is full.
      The constant invertible map $A_{\mathbf R}$ also preserves fullness, embeddedness, and descent to any quotient of $F$.

      If $A_{\mathbf R}F=O A_{\mathbf R'}F$ for orthogonal $O$, fullness implies $A_{\mathbf R}=O A_{\mathbf R'}$.
      Hence $A_{\mathbf R}^*A_{\mathbf R}=A_{\mathbf R'}^*A_{\mathbf R'}$, so $R_i=R_i'$ for every $i$.
      Finally, $\sum_i t_iR_i=0$ leaves $(q-1)\dim\mathcal K_0(f)$ free parameters, and the condition $I+R_i>0$ holds in a neighborhood of zero.
\end{proof}

\begin{example}[A ten-parameter family with one metric]
      Take $m=3$, $d_1=d_2=2$, $k_1=k_2=2$ in Proposition \ref{prop:har}.
      Then $\lambda_2=8$, $\dim E_2=9$, and $t_1=t_2=\frac12$.
      We obtain a full minimal embedding
$$
          F([x],u,v)=2^{-\frac12}
          \bigl(\phi_2(x)\otimes u,\phi_2(x)\otimes v\bigr)
          :\mathbb RP^3\times\Sn^2\times\Sn^2\longrightarrow\Sn^{53}
$$
      with $g=\frac83g_{\mathbb RP^3}+\frac12g_{\Sn^2}+\frac12g_{\Sn^2}$,
          where $g_{\mathbb RP^3}$ is induced from the unit round sphere.
      The identity $(x_1x_2)(x_3x_4)=(x_1x_3)(x_2x_4)$ gives a nonzero element of $\mathcal K_0(\phi_2)$.
      By \cite[Theorem~A]{Toth2002}, multiplication maps $\operatorname{Sym}^2E_2$ onto $E_0\oplus E_2\oplus E_4$ with kernel
          $\mathcal K_0(\phi_2)$, so
          $\dim\mathcal K_0(\phi_2)=\frac{9\cdot10}{2}-(1+9+25)=10$.
      Taking $(R_1,R_2)=(R,-R)$ for sufficiently small $R\in\mathcal K_0(\phi_2)$,
          we obtain a ten-parameter family of minimal embeddings of this fixed metric into $\Sn^{53}$.
      For a fixed source parametrization, these embeddings are pairwise inequivalent under orthogonal transformations of the target.
\end{example}

\Needspace{9\baselineskip}
\section{Clifford families and their first eigenfunctions}\label{sec:cli}

We combine the position and unit normal of a Clifford focal manifold with a complete harmonic map applied to the normal coefficients.
      The resulting minimal embeddings share one fixed unit normal bundle as their source.
      For sufficiently large harmonic degree, the first eigenspace of the induced metric determines a fibration of this bundle.
      We use this to describe its isometries and distinguish the metrics.
      The final example also gives intrinsically nonhomogeneous sources.

\subsection{Construction from a Clifford system}

      Let $\{P_0,\ldots,P_m\}$ be a symmetric Clifford system on $\mathbb R^{2\ell}$,
          so $P_\alpha P_\beta+P_\beta P_\alpha=2\delta_{\alpha\beta}I$, and write $P_u=\sum_\alpha u_\alpha P_\alpha$.
      We assume $m\geq1$ and $m_2=\ell-m-1>0$.
      The associated focal manifold
$$
          B=M_+=\{x\in\Sn^{2\ell-1}:\langle P_\alpha x,x\rangle=0,\ 0\leq\alpha\leq m\}
$$
      is closed, connected and minimal, with induced metric $g_B$ and dimension $d=2\ell-m-2=m+2m_2$; see \cite[\S\S~2.8, 4.2]{FKM}.
      The normal frame $(P_0x,\ldots,P_mx)$ identifies its unit normal bundle with $N=B\times\Sn^m$.
      On $N^n$, where $n=d+m=2\ell-2$, we set
\begin{equation}\label{eq:clifford-maps}
          f(x,u)=x,\qquad g(x,u)=P_ux,\qquad
          \pi(x,u)=(\pi_0(x,u),\ldots,\pi_m(x,u))=u.
\end{equation}
      Thus $f$ and $g$ record the position and the chosen unit normal.
      The projection $\pi:N\to\Sn^m$ records the normal coefficients and has fibers $B\times\{u\}$.

\begin{lemma}[The position and normal maps]\label{lem:cga}
      The map $H=\frac1{\sqrt2}(f,g):N\to\Sn^{4\ell-1}$ is a linearly full minimal embedding.
      For $V\in T_u\Sn^m$, put $L_uV=P_uP_Vx$.
      Its induced metric is
\begin{equation}\label{eq:cga}
          h_0((X,V),(X,V))
          =\left|X+\frac12L_uV\right|^2+\frac14|V|^2.
\end{equation}
      Each fiber is $(B,g_B)$ and is minimally embedded by $H$.
      The horizontal distribution of $\pi$ is
$$
          \mathcal H_{(x,u)}
          =\left\{\left(-\frac12L_uV,V\right):V\in T_u\Sn^m\right\}.
$$
\end{lemma}

\begin{proof}
      We have $dH(X,V)=\frac1{\sqrt2}(X,P_uX+P_Vx)$.
      The Clifford relations imply $L_uV\in T_xB$ and $|L_uV|=|V|$, from which we obtain \eqref{eq:cga} and $\mathcal H$.
      For fixed $u$, the fiber restriction of $H$ is a linear isometric inclusion of $B$ and is therefore minimal.

      Choose a geodesic frame $v_i$ at $u$ and put $K_i(x)=P_uP_{v_i}x$.
      These are tangent Killing fields on $B$ by \cite[\S~6.2]{FKM}.
      The volume density of $h_0$ is $2^{-m}$ times the product density, and $E_i=(-K_i,2v_i)$ form a horizontal orthonormal frame.
      At the chosen point, $\operatorname{div}_BK_i=0$ and $\partial_{v_i}(P_uP_{v_i})=P_{v_i}^2-P_u^2=0$.
      Hence $\operatorname{div}_{h_0}E_i=0$, and we obtain
$$
          \Delta_{h_0}f=\Delta_Bf+\sum_iK_i^2f=-(d+m)f=-nf,
$$
      since $(P_uP_{v_i})^2=-I$.
      Conjugation by $P_u$ preserves the Clifford span and therefore $B$.
      The involution $I(x,u)=(P_ux,u)$ exchanges $f,g$ and is an $h_0$-isometry.
      It follows that $\Delta_{h_0}g=-ng$, proving minimality.

      To prove that $B$ is linearly full, suppose a nonzero constant vector $z$ is orthogonal to its span.
      Then $z=P_{v(x)}x$ with $|v(x)|=|z|$, so $x=|z|^{-2}P_{v(x)}z$ would place the $d$-manifold $B$ in the image of an $m$-sphere, contrary to $d>m$.
      Comparing $u$ and $-u$ now proves fullness of $H$.
      Its first block determines $x$, and $|P_{u-u'}x|=|u-u'|$ determines $u$.
      Thus $H$ is an embedding.
\end{proof}

      To combine $H$ with an eigenmap of $\Sn^m$, we use the two-block balance from the ordinary Lawson--Osserman construction in \S\ref{sec:mmp}.
      We first carry out the calculation for a general minimal fibration.

\begin{proposition}[Matching over a minimal fibration]\label{prop:fib}
      Let $F_0:(E^n,g_0)\to\Sn^a$ be a linearly full minimal isometric embedding of a closed connected manifold.
      Let $\pi:(E,g_0)\to(C^m,\kappa g_C)$ be a surjective Riemannian submersion with connected fibers,
          where $0<m<n$ and $\kappa>0$.
      Assume that $F_0$ restricts to a minimal immersion on each fiber.
      Write $d=n-m$, $g_H=\pi^*g_C$ and $g_V=g_0-\kappa g_H$.
      Let $\psi_\lambda:C\to\Sn^b$ be a linearly full map satisfying
          $\psi_\lambda^*g_{\Sn^b}=\frac\lambda m g_C$ and $\Delta_{g_C}\psi_\lambda=-\lambda\psi_\lambda$, with $\lambda>0$.
      For $0<t<1$, set $F_t=(\sqrt t\,F_0,\sqrt{1-t}\,\psi_\lambda\circ\pi):E\to\Sn^{a+b+1}$.
      This embedding is minimal exactly when $\lambda>n\kappa$ and $t=t_\lambda:=\frac{d\lambda}{n(\lambda-m\kappa)}$.
      It is then linearly full and induces $h_\lambda=t_\lambda g_V+\frac\lambda n g_H$.
\end{proposition}

\begin{proof}
      Put $B_t=\kappa t+\frac\lambda m(1-t)$, so the induced metric is $h_t=tg_V+B_tg_H$.
      Since each fiber is minimal in the sphere, its mean curvature in $(E,g_0)$ also vanishes.
      Using both minimality assumptions, we compute the Laplacian after constant rescaling:
$$
          \Delta_{h_t}F_0=-\left(\frac dt+\frac{m\kappa}{B_t}\right)F_0,
          \qquad
          \Delta_{h_t}(\psi_\lambda\circ\pi)=-\frac\lambda{B_t}\psi_\lambda\circ\pi.
$$
      Equating both eigenvalues to $n$, we obtain $B_t=\frac\lambda n$ and $t=t_\lambda$.
      The condition $0<t_\lambda<1$ is equivalent to $\lambda>n\kappa$.
The first block makes $F_t$ an embedding.
      Since $\Delta_VF_0=-dF_0$, the mean of $F_0$ on each fiber is zero.
      Averaging a linear relation over each fiber separates the two full blocks and proves fullness.
\end{proof}

      The cubic construction in the author's unpublished draft \cite[\S~4]{Zhang2021} follows directly from the proposition.
      For the Hopf submersion $\pi:\Sn^3\to\Sn^2(\frac12)$, write $g_{\Sn^3}=g_H+\eta^2$.
      The classical cubic map
$$
          P_3(z,w)=(z^3,\sqrt3\,z^2w,\sqrt3\,zw^2,w^3)\in\Sn^7
$$
      induces $3g_H+9\eta^2$; see \cite[Example~3.4]{CDVV1996} for this cubic example.
      It descends to a full minimal embedding of $L(3,1)=\Sn^3/\mu_3$,
          where $\mu_3$ acts by scalar cube roots of unity.
We also write $\pi$ for the descended Hopf projection, whose fibers are mapped to great circles by $P_3$.
      Let $\psi_k:\Sn^2(\frac12)\to\Sn^{2k}$ be the complete degree-$k$ harmonic map and put $L_k=2k(k+1)$.
      Taking $F_0=P_3$, $\kappa=3$ and $\lambda=2L_k$ in Proposition \ref{prop:fib}, we obtain
$$
          \Phi_k=(a_k\psi_k\circ\pi,b_kP_3):L(3,1)\to\Sn^{2k+8},\qquad
          a_k^2=\frac{2L_k-9}{3(L_k-3)},\quad b_k^2=\frac{L_k}{3(L_k-3)},
$$
      a full minimal embedding for every $k\geq2$.

      Recall the maps $f$, $g$ and $\pi$ in \eqref{eq:clifford-maps}.
      For $k\geq1$, let $\phi_k:\Sn^m\to\Sn^{D_{m,k}-1}$ be the normalized complete harmonic map in \eqref{eq:har},
          with $\lambda_k=k(k+m-1)$ and $D_{m,k}=\binom{m+k}{k}-\binom{m+k-2}{k-2}$, taking the second binomial as zero for $k=1$.
      We use $\phi_k\circ\pi$ as the third block on $N$ and combine it with $H=\frac1{\sqrt2}(f,g)$ by Proposition \ref{prop:fib}.

\begin{proposition}[The matched Clifford family]\label{prop:clifford-family}
      The three blocks $f$, $g$ and $\phi_k\circ\pi$ admit a positive matched triple, with no auxiliary factors, exactly when $4\lambda_k>n$.
      The unique triple is $(s_k,s_k,c_k)$, where
$$
          t_k=\frac{4d\lambda_k}{n(4\lambda_k-m)},\qquad
          s_k=\frac{t_k}{2},\qquad
          c_k=1-t_k=\frac{m(4\lambda_k-n)}{n(4\lambda_k-m)}.
$$
      The resulting map
$$
          \mathcal G_k=\bigl(\sqrt{s_k}f,\sqrt{s_k}g,\sqrt{c_k}\,\phi_k\circ\pi\bigr)
          :N^n\longrightarrow\Sn^{4\ell+D_{m,k}-1}
$$
      is a linearly full minimal embedding with induced metric
\begin{equation}\label{eq:cmet}
          h_k((X,V),(X,V))
          =t_k\left|X+\frac12L_uV\right|^2+\frac{\lambda_k}{n}|V|^2.
\end{equation}
\end{proposition}

\begin{proof}
Lemma \ref{lem:cga} gives a Riemannian submersion $\pi:(N,h_0)\to(\Sn^m,\frac14g_{\Sn^m})$ with minimal fiber restrictions of $H$.
      By \eqref{eq:cga}, its vertical metric is $g_V((X,V),(X,V))=|X+\frac12L_uV|^2$.
      We apply Proposition \ref{prop:fib} with $F_0=H$, $\kappa=\frac14$ and $\lambda=\lambda_k$.
For uniqueness, the involution $I$ from Lemma \ref{lem:cga} exchanges $f,g$ and fixes $\pi$.
      It therefore sends any matched triple $(a,b,c)$ to $(b,a,c)$.
      Theorem \ref{thm:mat} implies that their induced metrics coincide, hence $(a-b)(f^*g_{\Sn}-g^*g_{\Sn})=0$.
      Evaluating on $(0,V)$, we obtain $(a-b)|V|^2=0$, hence $a=b$.
We combine the first two blocks as $H$ and obtain uniqueness from Proposition \ref{prop:fib}.
\end{proof}

\subsection{The first eigenspace and the fibration}

      We show that, for sufficiently large $k$, the first eigenspace of $h_k$ determines the fibers of $\pi$.
      This allows us to determine the isometry group and distinguish the metrics for different degrees.
      For the metrics in Proposition \ref{prop:clifford-family}, functions pulled back from $\Sn^m$ satisfy
$$
          \Delta_{h_k}(\varphi\circ\pi)=\frac n{\lambda_k}(\Delta_{\Sn^m}\varphi)\circ\pi,
          \qquad \Delta_{h_k}\pi=-\frac{nm}{\lambda_k}\pi.
$$
      Put $b_k=\frac{\lambda_k}{n}$ and $\nu_1=\lambda_1(B,g_B)$. We use the threshold
$$
          \Lambda_B=\max\left\{\frac n4,\frac m4+\frac{md}{\nu_1}\right\}.
$$

\begin{theorem}[Intrinsic determination of the fibration]\label{thm:crec}
      If $\lambda_k>\Lambda_B$, then
\begin{equation}\label{eq:cfst}
          \lambda_1(h_k)=\frac{nm}{\lambda_k},\qquad
          E_{\lambda_1(h_k)}(h_k)=\operatorname{span}\{\pi_0,\ldots,\pi_m\}.
\end{equation}
      The fibers of $\pi$ are exactly the common level sets of the first eigenfunctions.
      The full isometry group consists of the diffeomorphisms $\Phi$ satisfying
          $\pi\circ\Phi=A\circ\pi$ for some $A\in\operatorname O(m+1)$,
          preserving $\mathcal H$, and restricting to $g_B$-isometries between fibers.
      In particular, this group is independent of $k$ as a subgroup of $\operatorname{Diff}(N)$ throughout this range.
\end{theorem}

\begin{proof}
      The volume density of $h_k$ is $t_k^{\frac{d}{2}}b_k^{\frac{m}{2}}$ times the product density.
      Thus the $L^2$-orthogonal complement of base pullbacks consists of functions whose mean on each fiber is zero.
      Both spaces are invariant under the Laplacian.
      On the space of base pullbacks, the least positive eigenvalue is $\frac m{b_k}$.
      On its complement, the fiberwise Poincar\'e inequality bounds the Rayleigh quotient below by $\frac{\nu_1}{t_k}$.
      The hypothesis ensures positive matching coefficients and $\frac m{b_k}<\frac{\nu_1}{t_k}$, proving \eqref{eq:cfst}.

      The functions $\pi_i$ are $L^2$-orthogonal with equal norm and satisfy $\sum_i\pi_i^2=1$.
      For any orthonormal basis $\varphi=(\varphi_0,\ldots,\varphi_m)$ of the first eigenspace, we therefore have
          $\frac\varphi{|\varphi|}=A\pi$ for some $A\in\operatorname O(m+1)$.
      Hence every isometry preserves the fibers and their orthogonal complement $\mathcal H$, and restricts to a $g_B$-isometry between fibers.
      Conversely, the stated diffeomorphisms preserve both $g_V$ and $\pi^*g_{\Sn^m}$, and therefore every $h_k$.
\end{proof}

      The same fiber decomposition distinguishes the metrics arising from different degrees.

\begin{corollary}[Distinct metrics]
      The metrics $h_k$ with $\lambda_k>\Lambda_B$ are pairwise nonhomothetic.
      If $m\geq2$ and $m_2>1$, they are pairwise not locally homothetic.
\end{corollary}

\begin{proof}
      Suppose $\Phi^*h_j=c h_k$.
      Pullback by $\Phi$ identifies the first eigenspaces, so $\pi\circ\Phi=T\pi$ for a linear map $T$.
      Since $\pi$ is onto the unit sphere, $T\in\operatorname O(m+1)$.
      On fibers, we have $\Phi^*(t_jg_B)=ct_kg_B$, and comparison of volumes yields $t_j=ct_k$.
      On their orthogonal complements, we obtain $b_j=cb_k$.
      Thus $\frac{b_j}{t_j}=\frac{b_k}{t_k}$, and the formula $\frac{b_k}{t_k}=\frac{4\lambda_k-m}{4d}$ forces $j=k$.

      For $0<\theta<\frac\pi4$, the normal tube map $(x,u)\mapsto\cos\theta\,x+\sin\theta\,P_ux$
          identifies $N$ with a regular isoparametric hypersurface.
      By \cite[\S~4.2, p.~485]{FKM}, this hypersurface is an $\Sn^{m_2}$-bundle over $M_-$,
          and $M_-$ is an $\Sn^{\ell-1}$-bundle over $\Sn^m$.
      When $m\geq2$ and $m_2>1$, the homotopy exact sequences imply that $N$ is simply connected.
      Since the metrics are complete and real analytic, a local homothety extends, after rescaling, to a global isometry
          by \cite[Vol.~I, Cor.~6.4, p.~256]{KobayashiNomizuI}.
      The preceding argument then proves the local assertion.
\end{proof}

      The restriction $m\geq2$ is necessary for the local assertion.
      When $m=1$, write $u=(\cos\theta,\sin\theta)$ and $J=P_0P_1$.
      The change $y=e^{\frac{\theta J}{2}}x$ puts \eqref{eq:cmet} locally in the form $t_kg_B(y)+\frac{\lambda_k}{n}d\theta^2$,
          so all admissible members are locally homothetic.

      The shape operators of $B$ in unit normal directions have norm at most one by \cite[\S~4.5(i)]{FKM}.
      The Gauss equation then gives $\operatorname{Ric}_B\geq2(m_2-1)g_B$.
      For $m_2>1$, Lichnerowicz's estimate \cite{Lichnerowicz} yields $\nu_1\geq\frac{2d(m_2-1)}{d-1}$.
      Consequently, it suffices to have
\begin{equation}\label{eq:cexp}
          \lambda_k>\max\left\{\frac n4,\frac m4+\frac{m(d-1)}{2(m_2-1)}\right\}.
\end{equation}

\begin{example}[A nonhomogeneous Clifford family]
      For multiplicities $(m,m_2)=(3,8)$, we have $\ell=12$, $d=19$, $n=22$ and $D_{3,k}=(k+1)^2$.
      Positive matching is possible exactly for $k\geq2$.
      Each such degree satisfies \eqref{eq:cexp}, since
          $\lambda_k=k(k+2)\geq8>\max\{\frac{11}{2},\frac{129}{28}\}$.
      Thus the full minimal embeddings $\mathcal G_k:N^{22}\to\Sn^{47+(k+1)^2}$ induce metrics that are pairwise not locally homothetic.
      Their first eigenspaces are $\operatorname{span}\{\pi_0,\ldots,\pi_3\}$, with eigenvalue $\frac{66}{k(k+2)}$.
      Together with Theorem \ref{thm:crec}, this proves Theorem C.

      These metrics are also intrinsically nonhomogeneous.
      By \cite[\S~5.8, p.~493]{FKM}, the common kernel of the shape operators $A_0,\ldots,A_3$ of $B$
          has minimum dimension zero and maximum dimension three.
      The Gauss equation gives $\operatorname{Ric}_B^{\sharp}=18I-\sum_{\alpha=0}^3A_\alpha^2$,
          so this kernel is the Ricci eigenspace with eigenvalue $18$.
      Its varying dimension shows that $(B,g_B)$ is not homogeneous.
      By Theorem \ref{thm:crec}, every isometry of $h_k$ preserves the fibration.
      If $(N,h_k)$ were homogeneous, the stabilizer of each fiber would act transitively on that fiber, contradicting the nonhomogeneity of $(B,g_B)$.
\end{example}

\appendix

\section{First-jet separation for spherical harmonics}
\label{app:jet}

      We prove Lemma \ref{lem:jet} by showing that the maps $f\otimes g\mapsto fg$
          and $f\otimes g\mapsto df\wedge dg$ together determine every tensor in $H_p\otimes H_q$.
      Let $\omega$ be the outward area form on $\Sn^2$, with Poisson bracket $df\wedge dg=\{f,g\}\omega$.
      Multiplication $\mu(f\otimes g)=fg$ and the bracket $\mathcal P(f\otimes g)=\{f,g\}$ are $\operatorname{SO}(3)$-equivariant.
      On each summand of the multiplicity-free decomposition
          $H_p^{\mathbb C}\otimes H_q^{\mathbb C}
          =\bigoplus_{\ell=|p-q|}^{p+q}H_\ell^{\mathbb C}$, we have
\begin{equation}\label{eq:sel}
\begin{aligned}
          \mu|_{H_\ell^{\mathbb C}}\ne0
          &\quad\Longleftrightarrow\quad p+q-\ell\text{ is even},\\
          \mathcal P|_{H_\ell^{\mathbb C}}\ne0
          &\quad\Longleftrightarrow\quad p+q-\ell\text{ is odd}.
\end{aligned}
\end{equation}
	      The antipodal map reverses $\omega$, so the two operations have parities $(-1)^{p+q}$ and $(-1)^{p+q+1}$, respectively.
      This proves the vanishing assertions.

	      To prove nonvanishing, we may assume $p\geq q$ and set $r=p-q$.
      In coordinates $x=(\sqrt{1-z^2}\cos\varphi,\sqrt{1-z^2}\sin\varphi,z)$, let $u=x_1+\mathrm i x_2$, $v=x_1-\mathrm i x_2$.
      Then $\omega=d\varphi\wedge dz$, $\{u,v\}=-2\mathrm i z$, and the harmonic polynomials $u^p,v^q$ satisfy
$$
\begin{aligned}
          u^pv^q&=e^{\mathrm i r\varphi}(1-z^2)^{\frac{r}{2}}(1-z^2)^q,\\
          \{u^p,v^q\}
          &=-2\mathrm i pq\,e^{\mathrm i r\varphi}(1-z^2)^{\frac{r}{2}}
            z(1-z^2)^{q-1}.
\end{aligned}
$$
      Up to scale, the weight-$r$ vector in $H_{r+j}^{\mathbb C}$ is $e^{\mathrm i r\varphi}(1-z^2)^{\frac{r}{2}}P_j(z)$,
          where $P_j$ is the monic Gegenbauer polynomial for $w_r(z)=(1-z^2)^r$.

      Separating variables in the spherical Laplacian, we obtain $L_rP_j=-j(j+2r+1)P_j$, with $L_r=(1-z^2)\partial_z^2-2(r+1)z\partial_z$.
      Since $w_rL_rf=((1-z^2)^{r+1}f')'$, integration by parts shows that $L_r$ is self-adjoint for $w_r\,dz$.
      The boundary terms vanish for polynomials.
      For $\varepsilon\in\{0,1\}$, define
          $F_a^\varepsilon=z^\varepsilon(1-z^2)^a$ and $c_a^\varepsilon=(2a+\varepsilon)(2a+\varepsilon+2r+1)$.
      Differentiating directly, we get
$$
          L_rF_a^\varepsilon
          =4a(a+r)F_{a-1}^\varepsilon-c_a^\varepsilon F_a^\varepsilon
          \qquad(a\geq1).
$$
      Fix $s\geq0$ and set $J_a=\int_{-1}^1w_rF_a^\varepsilon P_{2s+\varepsilon}\,dz$.
      Since $P_{2s+\varepsilon}$ is monic, $F_s^\varepsilon-(-1)^sP_{2s+\varepsilon}$ has lower degree.
      By orthogonality and self-adjointness, we obtain
$$
          J_s=(-1)^s\|P_{2s+\varepsilon}\|_{L^2(w_r)}^2\ne0,
          \qquad
          (c_a^\varepsilon-c_s^\varepsilon)J_a
          =4a(a+r)J_{a-1}\quad(a>s).
$$
      Since $c_a^\varepsilon-c_s^\varepsilon
          =2(a-s)(2a+2s+2\varepsilon+2r+1)>0$ for $a>s$,
          it follows inductively that $J_a\ne0$ for every $a\geq s$.
      For $(\varepsilon,a)=(0,q)$, the product has nonzero components in all degrees $r+2s$, $0\leq s\leq q$.
      For $(\varepsilon,a)=(1,q-1)$, the bracket has nonzero components in all degrees $r+2s+1$, $0\leq s<q$.
      This proves \eqref{eq:sel}.

      A nonzero equivariant map on an irreducible summand is injective, and components of distinct harmonic degrees are linearly independent.
      Thus $(\mu,\mathcal P)$ is injective on $H_p^{\mathbb C}\otimes H_q^{\mathbb C}$ and its real form.
      Choose real harmonic bases $(Y_\alpha)$, $(Z_\beta)$ orthonormal for probability area measure and,
          for the matrix $(t_{\beta\alpha})$ of $T$, define $A_T=\sum_{\alpha,\beta}t_{\beta\alpha}Y_\alpha\otimes Z_\beta$.
      The lemma's hypotheses are equivalent to
$$
          \mu(A_T)=0,\qquad
          \alpha_T:=\sum_{\alpha,\beta}t_{\beta\alpha}Z_\beta\,dY_\alpha=0.
$$
      Since $d\alpha_T=-\mathcal P(A_T)\omega$, injectivity implies $A_T=0$, and hence $T=0$. This completes the proof.

\begingroup\small
\begin{thebibliography}{99}
\setlength{\itemsep}{0pt}
\setlength{\parsep}{0pt}
\bibitem{Bryant}
R. L. Bryant, \emph{Minimal surfaces of constant curvature in $S^n$}, Trans. Amer. Math. Soc. \textbf{290} (1985), no. 1, 259--271.

\bibitem{Calabi1967}
E. Calabi,
\href{https://doi.org/10.4310/jdg/1214427884}{\emph{Minimal immersions of surfaces in Euclidean spheres}},
J. Differential Geom. \textbf{1} (1967), 111--125.

\bibitem{CDVV1996}
B.-Y. Chen, F. Dillen, L. Verstraelen and L. Vrancken,
\href{https://doi.org/10.1017/S0308210500030651}{\emph{An exotic totally real minimal immersion of $S^3$ in $\mathbb{CP}^3$ and its characterisation}},
Proc. Roy. Soc. Edinburgh Sect. A \textbf{126} (1996), no. 1, 153--165.

\bibitem{doCarmoWallach}
M. P. do Carmo and N. R. Wallach,
\href{https://doi.org/10.2307/1970752}{\emph{Minimal immersions of spheres into spheres}},
Ann. of Math. (2) \textbf{93} (1971), 43--62.

\bibitem{FKM}
D. Ferus, H. Karcher and H.-F. M\"unzner,
\href{https://doi.org/10.1007/BF01219082}{\emph{Cliffordalgebren und neue isoparametrische Hyperfl\"achen}},
Math. Z. \textbf{177} (1981), no. 4, 479--502.


\bibitem{KobayashiNomizuI}
S. Kobayashi and K. Nomizu,
\emph{Foundations of Differential Geometry}, Vol. I,
Interscience, New York, 1963.

\bibitem{Lawson1970}
H. B. Lawson, Jr.,
\emph{Complete minimal surfaces in $S^3$},
Ann. of Math. (2) \textbf{92} (1970), no. 3, 335--374.

\bibitem{LawsonOsserman}
H. B. Lawson, Jr. and R. Osserman,
\emph{Non-existence, non-uniqueness and irregularity of solutions to the minimal surface system},
Acta Math. \textbf{139} (1977), 1--17.


\bibitem{Lichnerowicz}
A. Lichnerowicz,
\emph{G\'eom\'etrie des groupes de transformations},
Travaux et Recherches Math\'ematiques, III,
Dunod, Paris, 1958.

\bibitem{LvWangXie2024}
Y. L\"u, P. Wang and Z. Xie,
\href{https://doi.org/10.1007/s00208-024-02799-8}{\emph{Classification of minimal immersions of conformally flat $3$-tori and $4$-tori into spheres by the first eigenfunctions}},
Math. Ann. \textbf{390} (2024), 2235--2280.

\bibitem{ParkUrakawa}
J.-S. Park and H. Urakawa,
\href{https://doi.org/10.1007/BF00147415}{\emph{Classification of harmonic mappings of constant energy density into spheres}},
Geom. Dedicata \textbf{37} (1991), no. 2, 211--226.

\bibitem{Takahashi}
T. Takahashi, \emph{Minimal immersions of Riemannian manifolds}, J. Math. Soc. Japan \textbf{18} (1966), 380--385.

\bibitem{Toth2002}
G. T\'oth,
\href{https://www.heldermann.de/JLT/JLT12/JLT122/jlt12034.htm}{\emph{Moduli for spherical maps and minimal immersions of homogeneous spaces}},
J. Lie Theory \textbf{12} (2002), no. 2, 551--570.

\bibitem{Toth2006}
G. T\'oth,
\href{https://doi.org/10.1090/S0002-9947-06-04081-5}{\emph{On the shape of the moduli of spherical minimal immersions}},
Trans. Amer. Math. Soc. \textbf{358} (2006), no. 6, 2425--2446.

\bibitem{XuYangZhang}
X. Xu, L. Yang and Y. Zhang,
\href{https://doi.org/10.1016/j.matpur.2019.01.019}{\emph{Dirichlet boundary values on Euclidean balls with infinitely many solutions for the minimal surface system}},
J. Math. Pures Appl. (9) \textbf{129} (2019), 266--300.

\bibitem{Zhang2021}
Y. Zhang,
\emph{Minimal warped product for minimal submanifolds in spheres},
unpublished draft, November 22, 2021.


\bibitem{ZhangMMP}
Y. Zhang,
\href{https://arxiv.org/abs/2609.24886}{\emph{Minimal Products via Weighted Spheres}},
arXiv:\allowbreak 2609.24886.
\end{thebibliography}
\endgroup
\end{document}